# Collatz's "3x+1" problem and iterative maps on interval

Wang Liang[1]

[1](Department of Control Science and Engineering, Huazhong University of Science and Technology, WuHan, 430074 , P.R.China)

**[Abstract]** In this paper, we convert Collatz map into a simple conjugate iterative maps defined in [0,1]. Such maps are more familiar to us and easier to deal with. Some new features of this map are observed by this method. An interesting heuristic proof is also presented. This research may provide a new venture to crack the "3x+1" problem.

**[Keywords]** "3x+1" problem，Symbolic dynamics

## 1 Introduction

A problem posed by L. Collatz in 1937, also called the $3x+1$ mapping, $3n+1$ problem, Hasse's algorithm, Kakutani's problem, Syracuse algorithm, Syracuse problem, etc:

$$T(x_n) = \begin{cases} \dfrac{x_{n-1}}{2}, & x_{n-1} \in even \\ \dfrac{3x_{n-1}+1}{2}, & x_{n-1} \in odd \end{cases}$$

"3x+1" Conjecture asserts that, starting from any positive integer n, repeat iteration of this function eventually produces the cycle(1,2). This "easy" and interesting problem has attracted many researchers, but till keeps in unsolved.

In my opinion, the main difficulty of this problem is that we can't find the proper method to study it. Now most well studied iterative maps like Logistic map ($x_{n+1} = 1 - ux_n^2, x \in [0,1]$) are all defined on intervals. A battery of tools have been developed to deal with them. But the "3x+1" map is defined in natural numbers, a discrete set. There seems still no efficient method to study such iterative maps.

So in this paper, we will convert "3x+1" maps into a conjugate map defined in an interval, which will have same dynamic as "3x+1" maps. Then we could apply some mature methods of dynamical system theory to study the "3x+1" maps.

## 2 Conjugacy function

**Definition 1:** Maps $f: X \to X$ and $g: Y \to Y$ are conjugate with conjugacy h if and only if there exists a bijection h such that commutes:

$$\begin{array}{ccc} X & \xrightarrow{f} & X \\ h \downarrow & & \downarrow h \\ Y & \xrightarrow{g} & Y \end{array}$$

Conjugacies preserve the dynamics of a map.

Corresponding author: E-mail: guoypm@hust.edu.cn

Here we wish to find a bijection $h: N \to [a,b]$ (N is natural number set) to convert the "3x+1" maps into a continued variable system. Because the N in "3x+1" map is divided into even and odd set, we also wish $h(x \in even)$ corresponds part of [a, b] and $h(x \in odd)$ is another part.

According to these principles, we construct h as follows:

$h: N \to [0,1]$

(1) For a natural number m, we write it in binary form: $\quad m = (a_n a_{n-1} \cdots a_1 a_0)_2, a_i \in \{0,1\}$

(2) Reverse the binary sequence (left/right): $\quad a_0 a_1 \cdots a_{n-1} a_n$

(3) Add a "0." Before sequence, turn it into a decimal fraction: $\quad 0.a_0 a_1 \cdots a_{n-1} a_n$

For example, select m=11, Its binary form is 1011. Then reverse it: 1101. Thirdly, add the "0.", it becomes $(0.1101)_2 = 0.8125$. Function h is symbols method and has no easy mathematic expression.

Obviously, for $m \in N, m < \infty$, $h(m)$ corresponds to the finite decimal fraction, which is only a small set of [0, 1]. For example, $\frac{1}{3} = (0.010101\cdots)_2$ has infinite length, so $h^{-1}\left(\frac{1}{3}\right) = \infty$. For irrational decimal fraction, it also corresponds to an infinite number. To construct a bijection, we should extend the domain of "3x+1" map from N to $N \cup \infty$.

So $h: N \cup \infty \to [0,1]$ is a bijection.

This extension has some difference with the 2-adic integers mentioned in other papers for "3x+1". For example, $(1111\cdots)_2 = -1$ in 2-adic, but here $(1111\cdots)_2 = \infty$, we just simply regard it as an infinite natural number. There is no difference for limited length number.

Moreover, even number must have the form $m = (a_n a_{n-1} \cdots a_1 0)_2$, $h(m) = (0.0 a_1 a_2 \cdots a_{n-1} a_n)_2$.

So $0 \le h(m \in even) < 0.5$. For odd number, $m = (a_n a_{n-1} \cdots a_1 1)_2, h(m) = (0.0 a_1 a_2 \cdots a_{n-1} a_n)_2$, so $0.5 \le h(m \in odd) < 1$

## 3 Conjugate maps of "3x+1"

After determining the h, we could calculate corresponding conjugate maps of "3x+1" as follows:

$$T(x) \xleftrightarrow{h} g(x) \Leftrightarrow T(x) = \begin{cases} T1(x) = \dfrac{x}{2}, x \in even \\ T2(x) = \dfrac{3x+1}{2}, x \in odd \end{cases} \xleftrightarrow{h} g(x) = \begin{cases} g1(x), x \in [0,0.5) \\ g2(x), x \in [0.5,1] \end{cases}$$

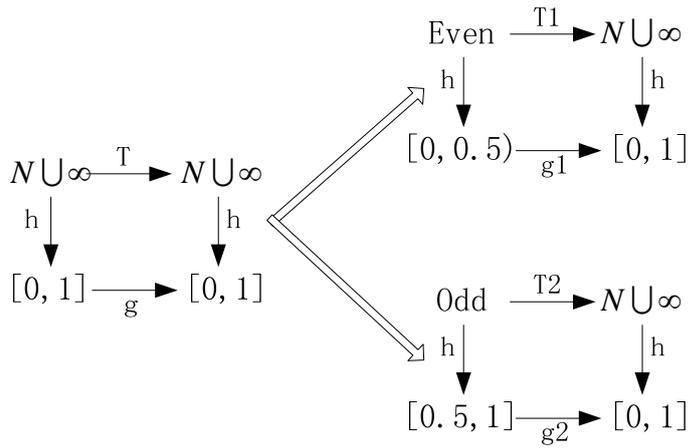

For $g1: h(x) \to h[T1(x)]$

Select $x = (a_n a_{n-1} \cdots a_2 a_1 0)_2 \in even$

$c1(x) = x/2 = (a_n a_{n-1} \cdots a_2 a_1 0)_2 / 2 = (a_n a_{n-1} \cdots a_2 a_1)_2$. For binary number, divided by 2 equal to right shift one bit.

$h(x) = (0.0 a_1 a_2 \cdots a_{n-1} a_n)_2$, $h(c1(x)) = (0.a_1 a_2 \cdots a_{n-1} a_n)_2$

We should calculate:

$?, g1: (0.0 a_1 a_2 \cdots a_{n-1} a_n)_2 \to (0.a_1 a_2 \cdots a_{n-1} a_n)_2$

Obviously, $2 \times (0.0 a_1 a_2 \cdots a_{n-1} a_n)_2 = (0.a_1 a_2 \cdots a_{n-1} a_n)_2$, So $g1(x) = 2x$

For $g2: h(x) \to h[T2(x)]$:

Select $x = (a_n a_{n-1} \cdots a_2 a_1 1)_2 \in odd$

$$T2(x) = \frac{3x+1}{2} = \frac{2x + x + 1}{2}$$

$$= \frac{2 \times a_n a_{n-1} \cdots a_3 a_2 a_1 1 + a_n a_{n-1} \cdots a_3 a_2 a_1 1 + 1}{2}$$

$$= a_n a_{n-1} \cdots a_3 a_2 a_1 1 + \frac{a_n a_{n-1} \cdots a_3 a_2 a_1 1 + 1}{2}$$

$$= a_n a_{n-1} \cdots a_3 a_2 a_1 1 + \frac{a_n a_{n-1} \cdots a_3 a_2 (a_1 + 1) 0}{2}$$

$$= a_n a_{n-1} \cdots a_3 a_2 a_1 1 + a_n a_{n-1} \cdots a_3 a_2 (a_1 + 1)$$
$$= a_n (a_{n-1} + a_n) \cdots (a_2 + a_3)(a_1 + a_2)(1 + a_1 + 1)$$
$$= a_n (a_{n-1} + a_n) \cdots (a_2 + a_3)(a_1 + a_2 + 1) a_1$$

(Attention: here the "+" is only a symbol. If we have to operate it, we should operate from left to right, but

not from right to left as normal. )

$h(x) = 0.1a_1 a_2 \cdots a_{n-1} a_n,$

$h(T2(x)) = 0.a_1(a_1 + a_2 + 1)(a_2 + a_3) \cdots (a_{n-1} + a_n)a_n$

We should find:

? $g2 : 0.1a_1 a_2 \cdots a_{n-1} a_n \to 0.a_1(a_1 + a_2 + 1)(a_2 + a_3) \cdots (a_{n-1} + a_n)a_n$

Unfortunately, it's very difficulty to give it a precise expression, but we could calculate some data of [h(x), h(c2(x)] and plot the rough figure of g2.

The figure of $g(x) = \begin{cases} g1(x) = 2x, x \in [0,0.5] \\ g2(x) = ?, x \in [0.5,1] \end{cases}$ is shown below:

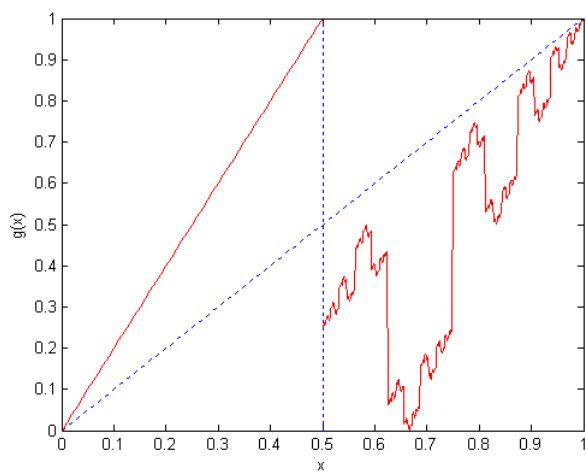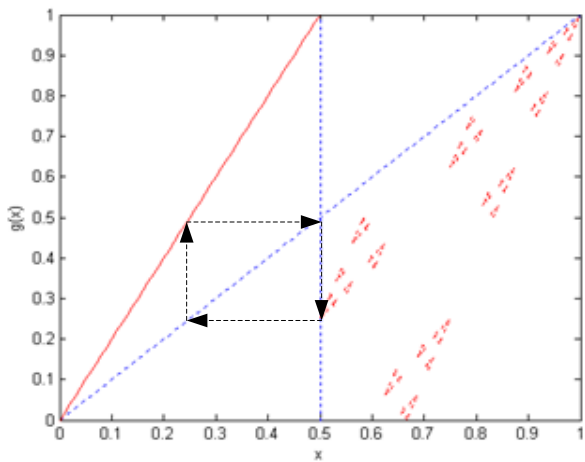

**Fig1.a** the figure of g(x). **Fig1.b** plot the g2(x) in dot form

## 4 A conventional proof method

Now the "3x+1" problem has been converted into the research of g(x), a map defined in an interval. According to the conjugacy theory, to prove the "3x+1" map is convergent to cycle (1,2), we just need to prove g(x) is convergent to cycle ( h(1), h(2) )=(0.5, 0.25).

Obviously, we only need consider the finite rational decimal in [0, 1] to solve "3x+1" problem. In this paragraph, we count all the numbers in.

Normally, if want to prove the orbits of a map $f$ are convergent, we should prove its fixed point or cycle is "attractive", which means they must satisfy:

For fixed point $x^*$: $\Delta = |f'(x^*)| < 1$.

For n-period cycle $x_1, x_2, \cdots, x_{n-1}, x_n$ : $\Delta = |\prod_{i=1}^{n} f'(x_i)| < 1$

Here we could also analyze g(x) following this way.

For any point in g1(x): $\frac{dg1}{dx} = 2$.

For any point in g2(x):

Select $m_1 = 0.1a_1 a_2 \cdots a_{n-1} a_n b_1 b_2 b_3 \cdots$ , $m_2 = 0.1a_1 a_2 \cdots a_{n-1} a_n c_1 c_2 c_3 \cdots$

When $n \to \infty$, $|m_1 - m_2| \to 0$. Assume $m_1 > m_2$, then $b_1 = 1, c_1 = 0$

$$\frac{dg2}{dx} = \frac{g2(m_1) - g2(m_2)}{m_1 - m_2} = \frac{g2(0.1a_1a_2 \cdots a_{n-1}a_n b_1 b_2 b_3 \cdots) - g2(0.1a_1a_2 \cdots a_{n-1}a_n c_1 c_2 c_3 \cdots)}{0.1a_1a_2 \cdots a_{n-1}a_n b_1 b_2 b_3 \cdots - 0.1a_1a_2 \cdots a_{n-1}a_n c_1 c_2 c_3 \cdots}$$

$$= \frac{0.a_1(a_1 + a_2 + 1) \cdots (a_{n-1} + a_n)(a_n + b_1)(b_1 + b_2) \cdots - 0.(a_1 + a_2 + 1) \cdots a_1(a_{n-1} + a_n)(a_n + c_1)(c_1 + c_2) \cdots}{0.000 \cdots 00(b_1 - c_1)(b_2 - c_2)(b_3 - c_3) \cdots}$$

$$= \frac{0.00 \cdots 00(a_n + b_1 - (a_n + c_1))(b_1 + b_2 - c_1 - c_2) \cdots}{0.000 \cdots 00(b_1 - c_1)(b_2 - c_2)(b_3 - c_3) \cdots}$$

$$= \frac{0.000 \cdots 0001(1 + b_2 - 0 - c_2) \cdots}{0.000 \cdots 00(1 - 0)(b_2 - c_2) \cdots} \geq 2$$

(Here we mark the n-th symbol as red color. In fact, we should prove whether derivative of $g(x)$ exist first. Here we just assume it exists. )

So $\Delta = |g'(x^*)| \geq 2 > 1$. There is no attractor in g(x).

For 2-period cycle, we have $\Delta = |g'(x_1)g'(x_2)| \geq 4 > 1$. That means there is also no attractive 2-period cycle.

In fact, if Collatz map is convergent to (1,2), the g(x) should also have an attractive 2-period cycle. This confused result shows such common method to judging the convergent or divergent of a map is not suitable for $g(x)$. We have to find another way to analyze it.

## 5 Special features of g(x)

Although g2(x) seems very complex, we could still find some special features of g(x), which may contain the new clews to crack it.

(1) $g2(x)$ is a fractal figure! The fractal is defined as "An irregular geometric object that is self-similar to its substructure at any level of refinement". This result is shown in Fig.2. It may be called Collatz fractal.

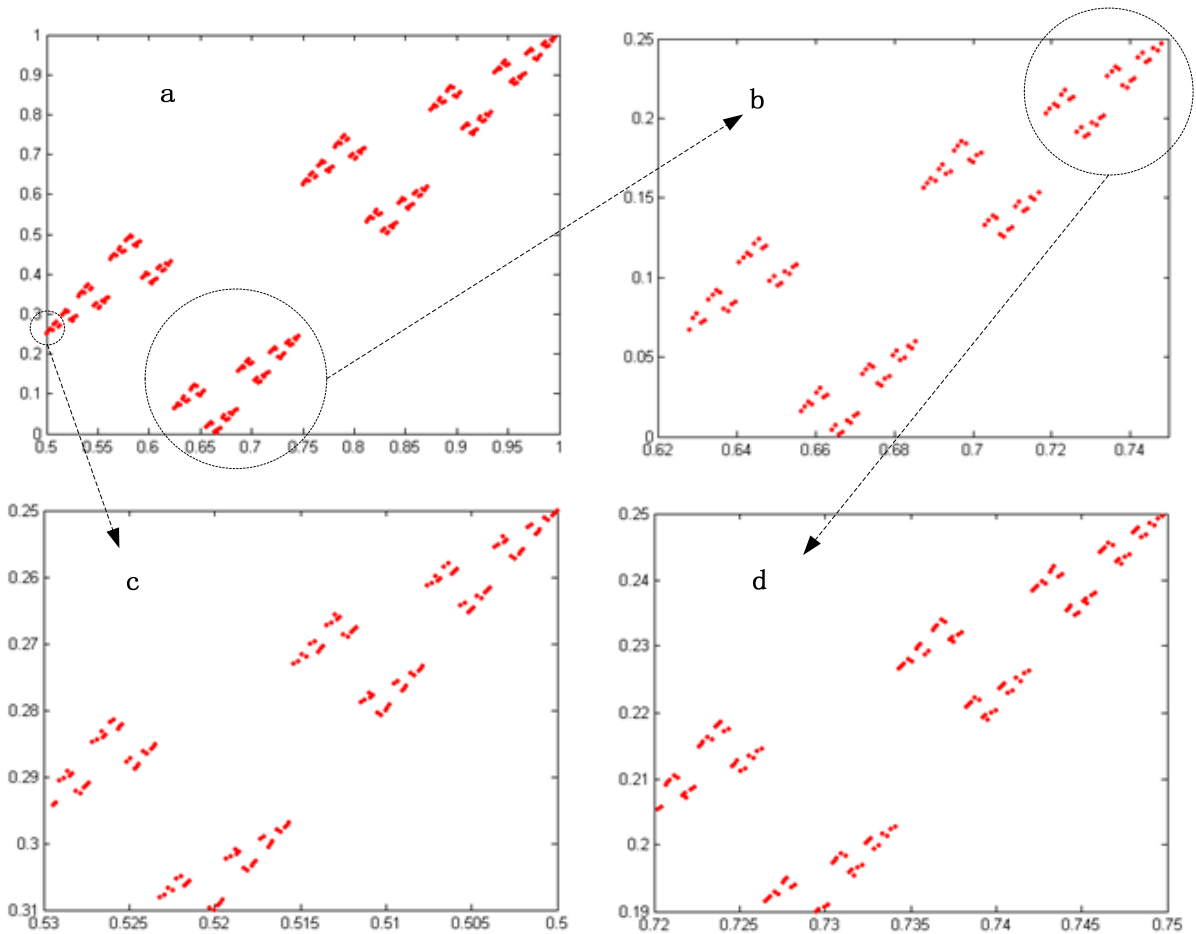

**Fig2.a** the figure of g2. **Fig2.b,c,d** part of g2. In Fig2.c, we reverse the X and Y axis.

(2) $g(x)$ and Bernoulli shift maps $B(x) = \begin{cases} 2x, 0 \leq x < \dfrac{1}{2} \\ 2x-1, \dfrac{1}{2} \leq x \leq 1 \end{cases}$ (Fig.3.a) are conjugate.

We know Bernoulli shift map play the kernel role in the research of continued variable iterative maps. So "3x+1" maps may be also a key to discrete variable iterative maps research. In fact, many complex number theory problems like some conjectures about prime distribution all could be regarded as discrete variable iterative problem.

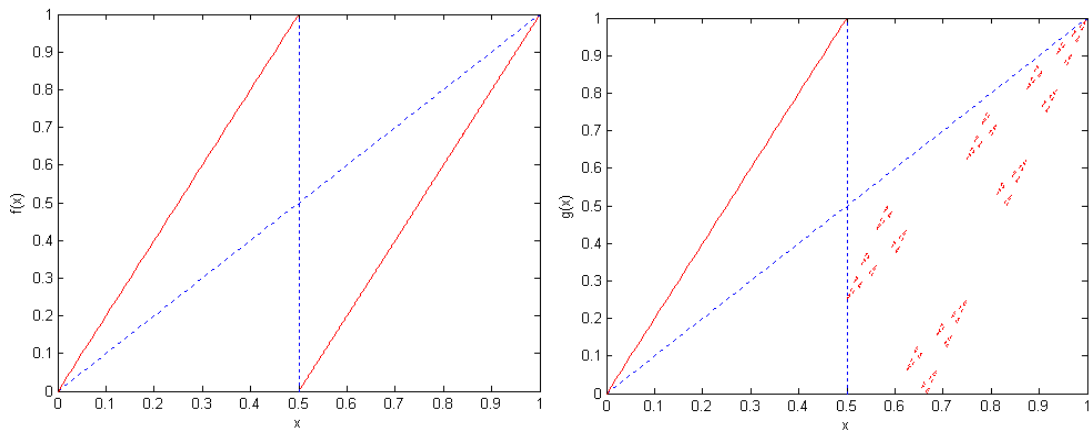

Fig.3 Bernoulli shift map and $g(x)$ are conjugate

This result could be proved as follows:

1. We can easily prove Bernoulli shift maps and another shift map $f(x) = \begin{cases} f1(x) = x/2, x \in even \\ f2(x) = (x-1)/2, x \in odd \end{cases}$ are conjugate with h : $f \xrightarrow{h} B$.

2. We have get $c \xrightarrow{h} g$ in this paper.

3. Lagarias has proved $\phi^{-1}$ (parity vector map) is a topological conjugacy between $c$ and $f$: $c \xrightarrow{\phi^{-1}} f$.

So ( $g \xrightarrow{h^{-1}} c, c \xrightarrow{\phi^{-1}} f, f \xrightarrow{h} B$ ) $\Rightarrow g \longrightarrow B$, $g$ and $B$ are conjugate.

(3) $g2(x)$ is discontinuous at every point in $[0.5,1]$.
Proof idea:
A example, $x_0 = 0.75 = (0.11)_2$ ,

$$g2: 0.1a_1a_2 \cdots a_{n-1}a_n \to 0.a_1(a_1+a_2+1)(a_2+a_3)\cdots(a_{n-1}+a_n)a_n$$

So $g2(0.11) = 0.1(1+1) = 0.101$
For $x^- < x_0$, $(x_0 - x^-) \to 0$, x could be written as $x^- = (0.101111\cdots 11)_2$, we have:
$g2(x^-) = g2(0.101111\cdots 11) = 0.0(0+1+1)(1+1)\cdots(1+1)1 = 0.00111\cdots 10 \to 0.01$
So $g2(x)$ is discontinues at point $x_0$.
In fact, for any point $x_0 = 0.a_0a_1a_2\cdots a_{n-1}1$, we could also write its close point as $x^- = 0.a_0a_1a_2\cdots a_{n-1}011\cdots 1$ and prove $g2(x^-) \neq g2(x_0)$. The strict proof could be given by the definition of continuous.

(4) is a bijection; $g2(x) \leq x$, (only when $x = (0.11\cdots 111)_2 = 1, g2(x) = x$) ; (0.5, 0.25) is a two period cycle of $g(x)$. The strict proof for these conclusions is not very difficult.

## 6 A heuristic proof

Here, we could give a heuristic proof for "3x+1" by basic dynamical system theory. This proof contains three steps:

**Theme 1:** Dividing the [0,1] into $2^n$ equal length parts, for any interval $I = [\frac{k}{2^n}, \frac{k+1}{2^n}], n \in N, k = 0,1,\cdots,2^n - 1$, we all have $|g(I)| = 2|I|$.

This theme means the g(x) double these intervals.

Main proof idea for this theme:

Obviously, $g1(x) = 2x$ could double any interval.

For $g2(x)$, this result could be shown by Fig.4.

We divide [0,1] into 4 equal parts in Fig1.a: The two independent parts of g2(x) are all in a rectangle with broad: long=1:2. It means $g2(x)$ also double the interval [0.5,0.75) and [0.75,1]. This conclusion is true for more detailed division (Fig.4.b, c, d).

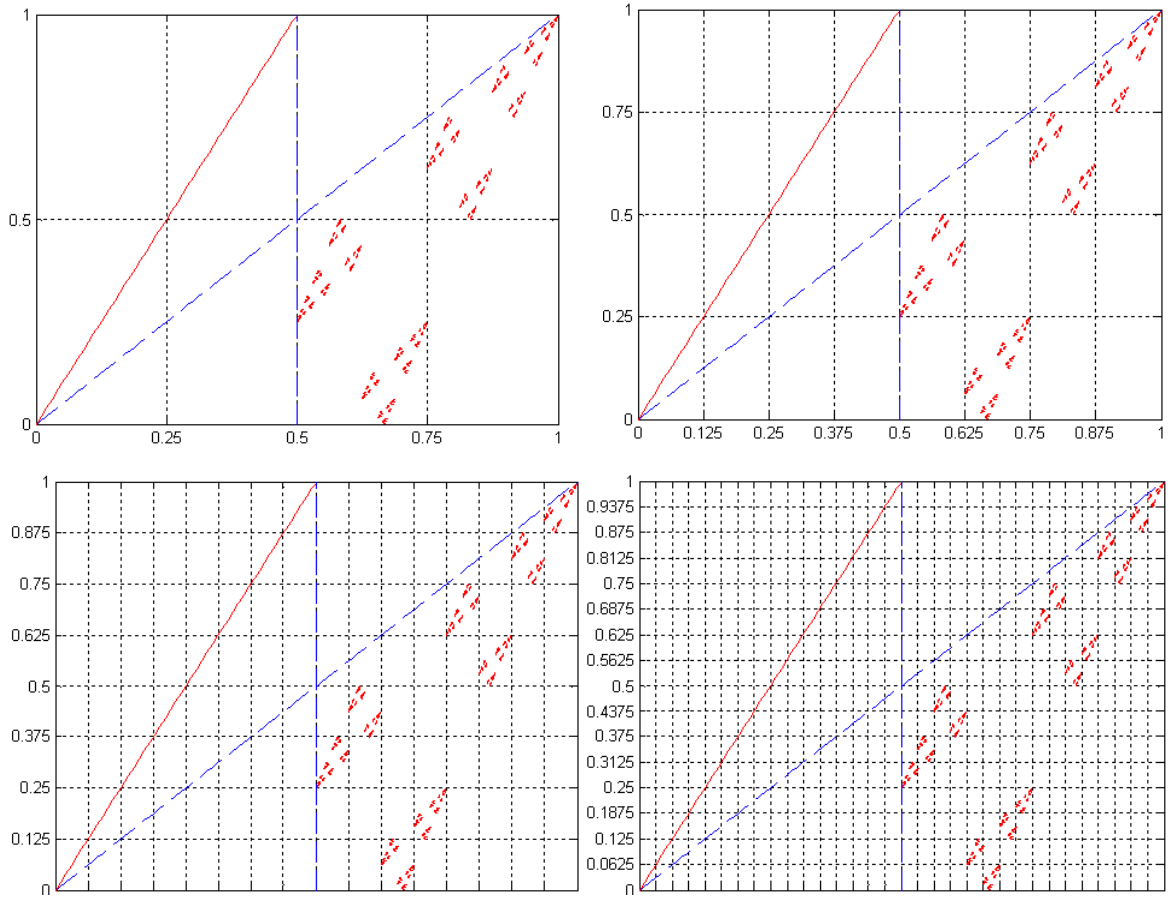

**Fig.4** $g(x)$ double the intervals

To give it a theory proof, we need find the max and min point of these intervals. We could find some basic patterns of interval iteration, which will be discussed in theme3.

Why do we select this division method?

For example, when we divide the [0,1] into $2^2 = 4$ equal parts, the endpoints of these intervals are $[0, (0.01)_2, (0.10)_2, (0.11)_2, 1]$. We could get all two length rational decimal. For $2^3 = 8$, the endpoints are $[0, (0.001)_2, (0.010)_2, \cdots, (0.110)_2, (0.111)_2, 1]$. So when dividing the [0,1] into more parts, we could get all the finite decimals in [0,1], which correspond to the natural numbers in "3x+1" maps. But for periodic circulating decimal and irrational decimal, we need divide the [0,1] into infinite parts $2^n$, n→∞.

Since it's very difficult to directly study the iteration of points, we could investigate the iteration of intervals. It's a common method in dynamical system research, which is also the main idea of our proof.

**Theme 2:** For any interval $I = \left[\frac{k}{2^n}, \frac{k+1}{2^n}\right], n \in N, k = 0,1,\cdots,2^n - 1$, we have $g^m(I) = [0,1], m < \infty$.

This theme means any interval $I$ could cover the whole $[0,1]$ after limited iterating. If theme1 is true, this conclusion will be obvious. We could also apply transfer matrix and transition graph/automata to describe this theme.

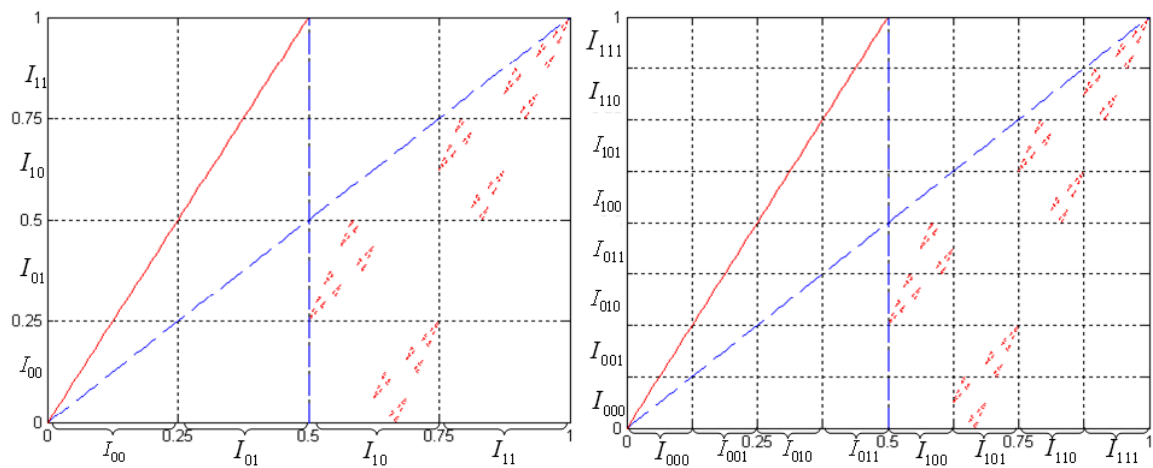

**Fig.5** Figure for transfer matrix

From Fig.5.a, we could get

$g(I_{00}) = I_{00} \cup I_{01}$
$g(I_{01}) = I_{10} \cup I_{11}$
$g(I_{10}) = I_{00} \cup I_{01}$
$g(I_{11}) = I_{10} \cup I_{11}$

For Fig.5.b, we could get the similar results.

Transfer matrix and transition graph is clearer method to describe these results, which is shown in Fig.6.

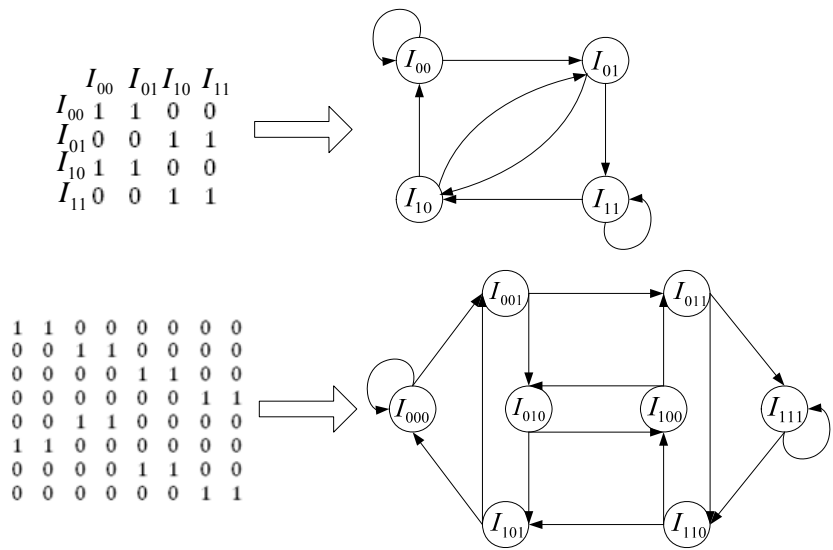

**Fig.6** transfer matrix and corresponding transition graph

From the transition graph, we could find that one interval could reach any other intervals through limited steps, which also means any interval could cover the whole $[0,1]$.

Normally, the transition graph could also be regarded as automata.

For example: $I_{00} \xrightarrow{g1} I_{01}$, $I_{11} \xrightarrow{g2} I_{11}$, etc. Here we remark g1 as '0' and g2 as '1'. So the transition graph could be converted into automata by marking their edges, which is shown in Fig.7.

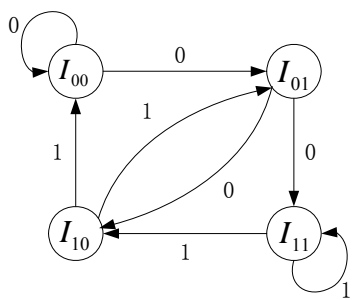

**Fig.7** Automata of g(x) maps

The relation of automata and "3x+1" is also discussed in some other papers. By our conjugate maps, this relation becomes more distinct.

In fact, the transfer matrix also features the iterate relation of residue classes of mod $2^n$ in "3x+1" map. For example, all the number in $I_{00}$ must begin with '0.00' (binary form). So according to conjugate relation, it corresponds to the natural numbers ending with '00' (binary form), which could be expressed as $x \equiv 0 \bmod 4$. By the similar method, we could also get:

$I_{01}: x \equiv 1 \bmod 4$
$I_{10}: x \equiv 2 \bmod 4$
$I_{11}: x \equiv 3 \bmod 4$

Related result could be shown by Fig.8.

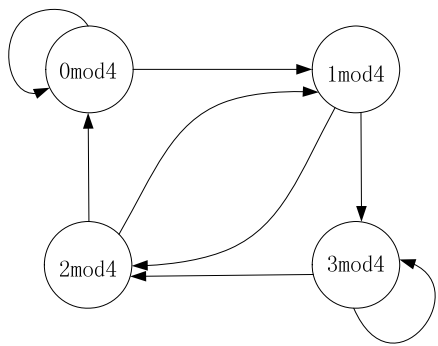

**Fig.8** dynamics of "3x+1" maps

**Theme 3:** For any interval $I = \left[\frac{k}{2^n}, \frac{k+1}{2^n}\right], n \in N, k = 0,1,\cdots,2^n - 1$, the end point of $I$ will reach the point 0.5 in x axis after limited iterating.

If $n < \infty$, the end points of $I$ is just the finite decimals in [0,1]. We also know (0.5, 0.25) is a two period cycle of $g(x)$. This theme means all finite decimals will be convergent to cycle (0.5, 0.25).

According to the conjugate theory, this theme also means all natural numbers will reach the cycle (1,2) after limited iterating of "3x+1" map. So Collatz's "3x+1" problem could be solved.

Main point for its proof:
(1) Three basic patterns of interval iterating(shown in Fig.9):

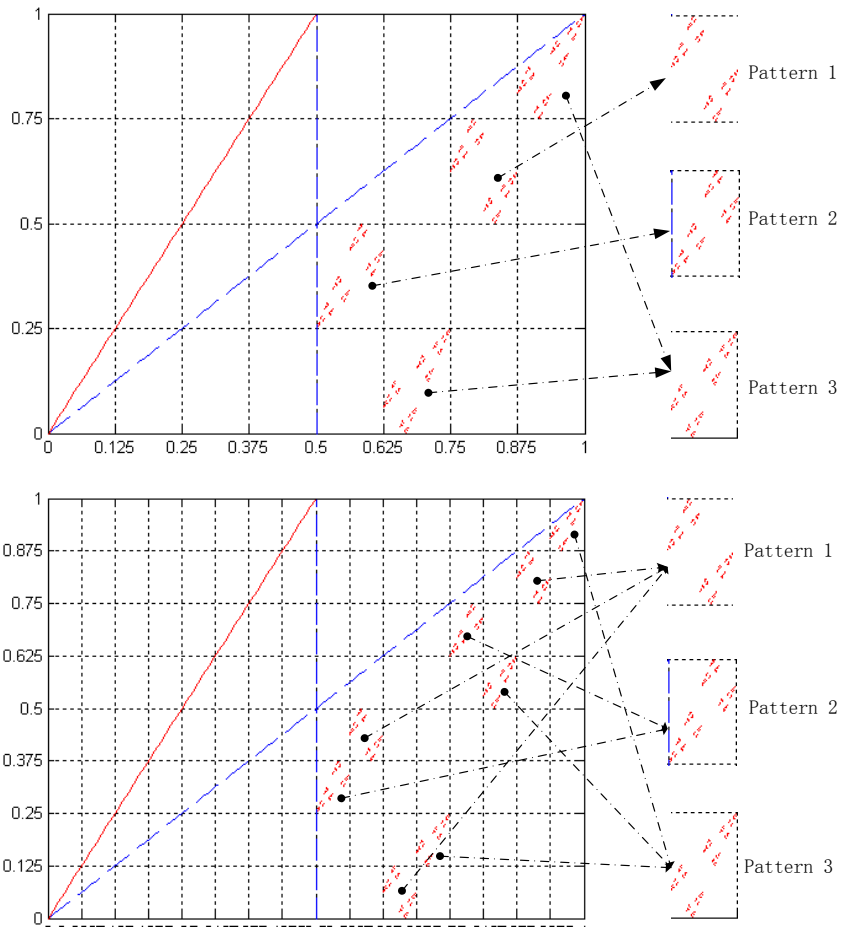

**Fig.9** basic patterns of interval iterating

In fact, these three patterns could also be induced into one pattern (pattern 2). It's shown in Fig.9:

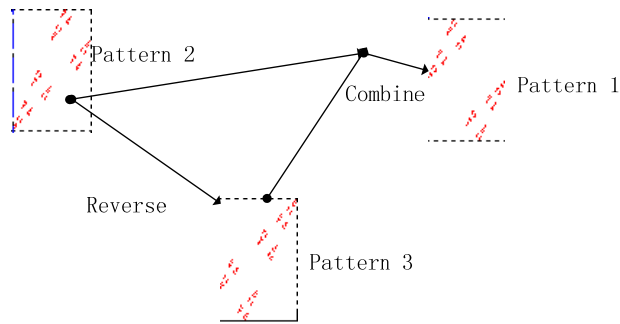

**Fig.9** relation of three basic patterns

So we could only use the figure of pattern 2 to construct the whole figure of $g2(x)$, which could contain any details. Such fractal structure may be the key for "3x+1" problem.

Furthermore, we only need pay attention to two features of pattern 2.

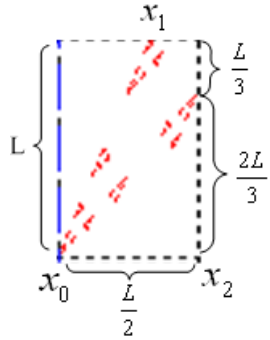

**Fig.9** basic pattern (pattern 2) of iteration

First, $g2(x_2) = ?$. $x_0$ is the minimal point of this interval. So if the $x_0$ is the endpoint of a $2^n$ division interval, it will reach another endpoint of one of these intervals after an iterating. For $x_2$, the experimental results show $|g2(x_2) - g2(x_0)| = \frac{2L}{3} = \frac{3}{2}|x_2 - x_0|$, which is shown in Fig.9. So if the $x_2$ is the endpoint of a $2^n$ division interval, it will reach the endpoint of a $2^{n+1}$ division interval after one iterating. This result could be observed in Fig.4.

Second, $x_1$ is the maximal value point of this interval. We want to know $x_1 = ?$. Experimental results shows $x_1$ may be not a finite decimal. The pattern of $x_1$ is the key for the proof of theme1.

(2) How to apply the induction method? Fractal structure of g2(x) implies we could use the induction proof method, but the key problem is how to construct the induction, by automata theory? Symbolic dynamics theory? Ergodic theory or something else?

# 7 Cousins problem of "3x+1"

Two simple cousins problem of "3x+1" are "5x+1" maps and "7x+1", which only change the c2(x)=(3x+1)/2 into c2(x)=(5x+1)/2 and c2(x)=(7x+1)/2. Some researches show they all contain the divergent orbits.

By the same method, we could also get their conjugate maps with function h. They are shown below:

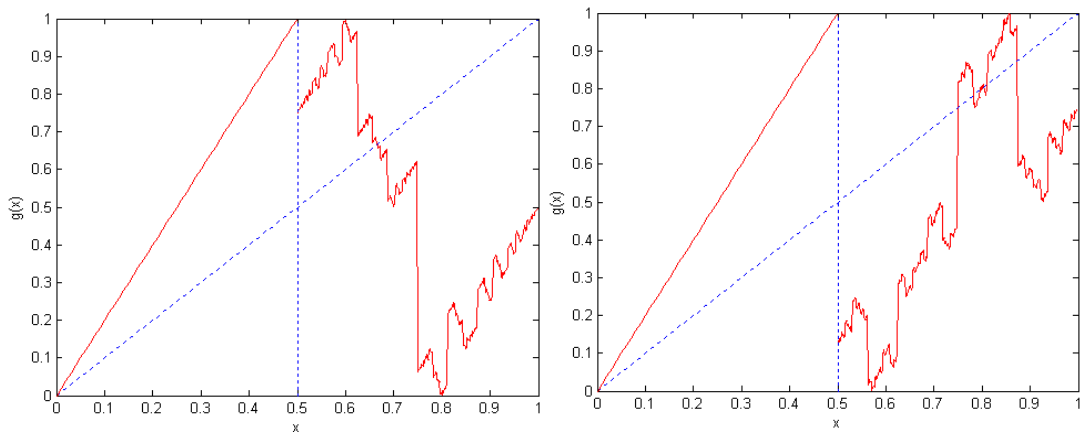

**Fig.10** the conjugate maps of "3x+1" are "5x+1" maps

We could find some parts of g2(x) of these two maps are above the y=x, but the g2(x) in "3x+1" is all below the y=x. This difference may have some relations with the convergent or divergent of these maps.

## 8 Conclusions

The discussion in this paper is not very strict. Our main aim is to discover the possibility to convert "3x+1" map into an iterative map on interval. Normally, such maps are easier to deal with. This research may open a new way for the solving of "3x+1" problem.

Here we construct the bijection between natural number and real number in [0,1] by a symbols method. We may also find a simpler bijection to analyze the "3x+1" problem in the future.

**Acknowledgement:** Thanks for the suggestion of Jeffrey Lagarias and Olivier Pirson.

# Appendix:

Two Matlab program for this paper:
**1 Function h:**
% must be saved as Collatz_h.m
% h:N --->[0,1]. Here n is a natural number
function h=Collatz_h(n)

% get the binary form of n
Bin_n=dec2bin(n);

%reverse the binary sequence of n
C_Bin_n=fliplr(Bin_n);

%turn it into a number in [0，1]
L=length(C_Bin_n);
Dec_g_n=0;
for i=1:L
    Dec_g_n=Dec_g_n+str2num(C_Bin_n(i))*2^(-i);
end
h=Dec_g_n;

---

**2 plot the figure of g(x)**
%Collatz_1.m
Clear
%plot the figure of g2(x)
n=1001;
A=[];
B=[];
for i=1:2:n
A=[A,Collatz_h(i)];    %A store the h(x)
f_n=(3*i+1)/2;        % it could be changed to f_n=(5*i+1)/2 , f_n=(7*i+1)/2 or other forms
B=[B,Collatz_h(f_n)]; %B store the h(c1(x))
end

%sort the data
[A_sort,index]=sort(A);
B_sort=B(index);

hold on
%figure of g2
plot(A_sort,B_sort,'r-')

```matlab
%figure of g1: y=2x
y=0:0.01:0.5;
plot(y,2*y,'r')

%figure of y=x
x=0:0.01:1;
plot(x,x,'b--')
```